%% file: derivform.tex
\pdfoutput = 1
\documentclass[12pt,a4paper,reqno]{amsart}

\usepackage{amssymb, amsmath, amsthm, amsxtra, amscd}

\input macro.tex

\def\droman{\,{\rm d}}
\def\dx{\droman x}
\def\dy{\droman y}
\def\dt{\droman t}

\title{On the antiderivative of inverse functions}

\begin{document}

\setbox0 = \hbox{\rm \large Micha\"el Bensimhoun,}
\setbox1 = \hbox{\rm Jerusalem, December 2013 \vrule height 15pt width 0pt }

\author{
\noindent\box0\\
\box1\\
}

%\bigskip
%\noindent {\sc 2000 AMS Subject Classification: } 47D06

\begin{abstract}
\small
One of the basics of calculus is the following proposition:
If $F$ and $G$ are antiderivatives 
of two real functions $f$ and $g$ resp.\@,
then an antiderivative of $\lambda f+g$ is $\lambda F + G$.
It may be surprising that such systematic an integration formula 
exists for the inverse of a function $f$, 
a fact that seems to have been discovered for the first time 
by Laisant in 1905 (\cite{La}), and seems to be still unsufficiently known.
More precisely, 
if $f$ is an invertible real function, 
and if $F$ is an antiderivative of $f$, then the antiderivative
of $f^{-1}$ is $xf^{-1}(x)-F\circ f^{-1}(x)+C$.
Laisant, and other authors after him (e.g.\@ \cite{Pa,St}), 
assumes that $f^{-1}$ is differentiable, in which case the proof 
of this formula is immediate.
Recently, it has been shown by Key that this additional
assumption is unnecessary (\cite{Ke}). 
In this paper, we give two different proofs of this result.  
The first proof, of geometrical spirit, relies to Fubini's theorem, while 
the second proof, purely analytic, is based on the Stieltjes integral. 
\end{abstract}
\medskip

\maketitle

\bigskip 

%================================================================== 
\parskip 3pt plus 1pt minus 1pt
\parindent 0pt  

Apparently, the aforementioned result has been first discovered by 
Laisant in 1905 (\cite{La}), but seems to be unconsciously used
by mathematicians, each time they need to compute the integral of 
the inverse of a function $f$.
This elementary theorem was rediscovered several time after
Laisant (see e.g.\@ \cite{Pa,St}), assuming furthermore that $f$ is differentiable, 
in which case the proof of the formula is immediate. 
It seems that Key was the first who proved the exactness of this
formula, without this additional assumption (\cite{Ke}). 
In his article, he uses it as an intermediate tool to prove that the shell and
disk methods for computing the volume of a solid of revolution
agree.  
 
The author of the present article, which was unaware of the papers of Laisant and Key, was 
lead to the same result by simple geometric considerations in 1999.
He then asked several mathematicians (some of renown) if this theorem
was not unknown to them. They replied that they are unaware of any previous
statement of this theorem, and that they are rather surprised that such elementary
a result is not included in any introducing book of calculus. 

Actually, if this theorem has already been discovered in the far past, 
it would be interesting
to understand why it is not widely known at the elementary level.
Also, it may be equally interesting to understand why this result was not sufficiently spread
since it was published in 1994.
The answer to this question may be that it is relatively recent (only 
twenty years). 
Another answer may be the way Key presented this theorem: 
The main topic, in his paper, is the proof of the agreement of the shell and 
disk methods. It is not entirely clear, from his article, 
that the intermediate proposition is a result comparable, in usefulness, 
to the integration by part formula.

The proof of Key is based on the very definition of the Riemann-Darboux integral,
\thatis the upper or lower limit of the Riemann-Darboux sums.
His argument is very clear and elegant, but it is based on the theory of 
Darboux sums, which is not always taught at the undergraduate level anymore.
In contrast, the first of the two alternative proofs offered in this article uses
more advanced tools of integral calculus (Fubini's theorem), 
but it does not assume any knowledge of the Riemann-Darboux integral,
which makes it independent of the kind of integral taught to students.
The second proof, essentially based on the change of variable theorem in 
the Stieltjes integral, has the 
advantage to be generalizable to more involved formulae.
Also, it should be observed that the argument needed to establish
the change of variable formula in the 
Stieltjes integral is very similar, in spirit, to the argument 
used in the proof of Key.  

Here is the precise statement of the aforementioned proposition.

\begin{theorem}
let $f: [a,b]\to [c,d]\incl \overline\bbR$  be a continuous 
and invertible function, with $a\geq -\infty$ and $b\leq +\infty$.
Then $f$ and $f^{-1}$ have antiderivatives, and if $F$ is an
antiderivative of $f$, then the antiderivative of $f^{-1}$ is
$$
G(x)= x f^{-1}(x)-F\circ f^{-1}(x)+C.
$$
\end{theorem}

\example
let $f(x)=e^x$; then by the formula above, one gets 
immediately 
$$
\int \log x \dx=x\log(x)-x+C.
$$
\medbreak

\begin{proof}
Since $f$ is continuous and invertible $[a,b]\to [c,d]$, 
it is an elementary result that $f^{-1}$ must be continuous. 
Therefore, the fundamental theorem of calculus implies that 
both $f$ and $f^{-1}$ have antiderivatives in $[a,b]$ and $[c,d]$ 
respectively.

Notice also that $f$ and $f^{-1}$ must be either both increasing, or 
both decreasing, as can be seen easily.
 
Because of the fundamental theorem of calculus, 
it suffices to show that
$$
G(y)=\int_{c}^y f^{-1}(t)\dt
$$
is equal to $y f^{-1}(y)-F\circ f^{-1}(y)+C$, 
with $C\in \bbR$ and 
$F(x)=\displaystyle\int_{f^{-1}(c)}^x f(x)\dx$.
Observe first that 
$$
f^{-1}(t) =\int_{f^{-1}(c)}^{f^{-1}(t)} \dx + f^{-1}(c),
$$
therefore 
$$
G(y) = \int_{c}^y \int_{f^{-1}(c)}^{f^{-1}(t)}\dx\dt + f^{-1}(c)y + C, 
\With C = -f^{-1}(c)c.
\eqno(1)
$$
In a system of axes $(t,x)$, let $D$ denote the domain 
between the horizontal straight $x=f^{-1}(c)$, the curve 
$x=f^{-1}(t)$, and the vertical straights $t=c$ and $t=y$. More formally,
if $f^{-1}$ is increasing,
$$
D = \set{(x,t)\in \bbR^2 \st c\leq t\leq y \andthat  
f^{-1}(c)\leq x \leq f^{-1}(t)},
$$
and if $f^{-1}$ is decreasing, 
$$
D = \set{(x,t)\in \bbR^2 \st c\leq t\leq y \andthat  
f^{-1}(t)\leq x \leq f^{-1}(c)}.
$$ 
Let us put $\epsi = 1$ if $f^{-1}$ is increasing, and $\epsi=-1$ if $f^{-1}$ is 
decreasing. Hence $\epsi=1$ if $f^{-1}(c)<f^{-1}(y)$ and $\epsi=-1$ 
if $f^{-1}(c)>f^{-1}(y)$.

By Fubini's theorem, 
$$
\int_{c}^y \int_{f^{-1}(c)}^{f^{-1}(y)}\dx\dt 
= \epsi \iint 1_D \dt\otimes{\rm d}x.  \eqno(2)
$$
For every fixed $x$ between $f^{-1}(c)$ and $f^{-1}(y)$, 
the function $\varphi(t)\mapsto 1_D(t,x)$ is 
equal to $1$ between $f(x)$ and $y$, and to $0$ otherwise. 
Indeed, assuming for example that $f$ is increasing (hence 
$f^{-1}$ is increasing), an element $(t,x)$ belongs to $D$
if and only if $c\leq t\leq y$ and $f^{-1}(c)\leq x\leq f^{-1}(t)$.  
Equivalently, $c\leq t\leq y$ and $c=f(a)\leq f(x)\leq t$. 
This can be put into the more compact form: $c\leq f(x)\leq t\leq y$.

Similarly, if $f$ is decreasing (hence so is $f^{-1}$), 
then $(t,x)$ belongs to $D$
if and only if $c\leq t\leq y$ and $f^{-1}(t)\leq x\leq f^{-1}(c)$.
Equivalently, $c\leq t\leq y$ and $c\leq f(x)\leq t$ (applying $f$ to the 
two sides of an inequality reverses the inequality since $f$ is decreasing). 
In compact form: 
$c\leq f(x)\leq t\leq y$. 
Thus, the above contention is established.
 
Fubini's theorem can be applied again:
\breaklines{
\int 1_D \dt\otimes{\rm d}x 
= \iint \varphi(t) \dt\dx 
=  \epsi\int_{f^{-1}(c)}^{f^{-1}(y)}\int_{f(x)}^y \dt \dx\\
= \epsi \int_{f^{-1}(c)}^{f^{-1}(y)} (y-f(x)) \dx
= \epsi \big[ y f^{-1}(y)-yf^{-1}(c) - F\big(f^{-1}(y)\big) + F\big(f^{-1}(c)\big) \big]\\
= \epsi \big[ yf^{-1}(y) - yf^{-1}(c) - F\big(f^{-1}(y)\big)\big] + C.
}

Taking eq.~(1) and~(2)  into account, it follows that
$$
G(y) = yf^{-1}(y) - F\big(f^{-1}(y)\big) + C, 
$$
as was to be shown.
\end{proof}
\bigbreak

\emph{Second proof:}
As above, we calculate 
$$
G(y)=\int_{c}^y f^{-1}(y)\dy.
$$
In what follows, the classic Stieltjes integral, as well as its famous
integration by part and change of variable theorems, will be used.

Since $f$ is an isomorphism $[a,b]\to [c,d]$, 
the change of variable formula is licit, that is, with  $y=f(x)$,
$$
G(y) = \int _{f^{-1}(c)}^x f^{-1}\big(f(x)\big)\droman f(x)
=  \int _{f^{-1}(c)}^x x \droman f(x).
$$
For the sake of clarity, let us put $g(x)=x$ and $\alpha=f^{-1}(c)$
(necessarily, $\alpha = a$ or $\alpha = b$).
The integration by part formula for the Stieltjes
integral is
$$
\int_{\alpha}^x g(x)\droman f(x)+\int_{\alpha}^x f(x)\droman g(x) 
= f(x)g(x)-f(\alpha)g(\alpha). 
$$
In other words, 
$$
\int_{\alpha}^x x \droman f(x) = f(x)x-\int_\alpha ^x f(x) \dx + C' 
\With C'=-f(\alpha)g(\alpha).
$$
Setting $F(x)=\int_\alpha^x f(x)\dx$, and substituting
$x=f^{-1}(y)$ inside this equation, there holds
$$
G(y) =  yf^{-1}(y) - F\big(f^{-1}(y)\big) + C, \With 
C = C'-F(\alpha).
$$ 
This ends the second proof of the theorem.

\remarks
\par
1) The same approach can be used to obtain more involved formulae 
containing inverse functions, such as 
$$
\int H\big(y, f^{-1}(y)\big) \dy 
= H\big(y, f^{-1}(y)\big)\,y - \int _\alpha^{f^{-1}(y)} f(x) \droman H(f(x),x) + C.
$$

2)
If it can be supposed that $f^{-1}$ is absolutely continuous, 
then a simpler proof can be given, based on known theorems about 
absolute continuity:
Indeed, $f^{-1}$ is invertible and monotonic by hypothesis, hence,
if $f^{-1}$ is also absolutely continuous, the function 
$$
G(x)=xf^{-1}(x)-F\circ f^{-1}(x)+C
$$
is also absolutely continuous,
because it is the sum, product, and composition of 
absolutely continuous and bounded functions. 
By a known theorem due to Lebesgue, $f^{-1}$ and $G$
are differentiable almost everywhere, and are the integral of their
derivative.
But at every point where $f^{-1}$ is differentiable, $G(x)$ is 
differentiable with derivative equal to $f^{-1}$.
This prove that the integral of $f^{-1}$ is precisely $G$.

%-----------------------------------------------------------------------------

% REFERENCES:

\bibliography{derivform}{}
\bibliographystyle{amsplain}
\end{document}

% Another possibility to include bibliography is:

\end{document}

%% file: macro.tex
\let\rm\normalfont

	\def\hbar{{\bar h}}

\def\bbR{\mathbb{R}}

\let\epsi=\varepsilon

\let\Phi=\phi
\let\phi=\varphi

\def\implies{\ifinner \Rightarrow \else \Longrightarrow \fi} 
\let\incl=\subseteq

\def\maps{\ifinner \mapsto \else \longmapsto \fi}

\def\andthat{\ \hbox{and}\ }

\def\st{\colon} % such that symbol

\def\enu#1 {\par \noindent (#1) }  

\def\set#1{\{#1\}}
\def\droman{{\mathrm d}} 

\def\dt{{\droman t}}

\def\dx{{\droman x}}
\def\dy{{\droman y}}

\newmuskip\halfthinmuskip
\newtheorem{theorem}{Theorem}[section]

\def\example{\bigskip {\noindent \bf Example: }}

\def\remarks{\bigskip {\noindent \bf Remarks: }}

\def\enum#1-#2{%
	\in\{\if\voidarg{#1}1,2\else #1\fi,\ldots,#2\}%
}

\protected\def\voidarg#1{%
	//\let\unlessb\empty\else\let\unlessb\unless\fi 
  \begingroup 
  	\aftergroup\expandafter\aftergroup\unlessb
  	\def\tempvar{#1}%
	\ifx\tempvar\empty \aftergroup \iftrueb 
	\else \aftergroup \iffalseb 
	\fi
  \endgroup
}  

\def\iffalseb{\iffalse}
\def\iftrueb{\iftrue}
\let\If=\if

\def\diffD{\mathop{\raise .05em\hbox{${\mathsurround 0pt\smallsetminus}$}}}%               
\def\diffS{\mathop{\raise .085em\hbox{$\mathsurround 0pt\scriptstyle\smallsetminus$}}}%

\def\restr{\raise -2 pt \hbox{$\m@th |$}}

\def\thatis{that~is, }
\def\With{\quad \mbox{\rm with}\quad}

\def\andthat{\ \mbox{\rm and}\ }

\def\breaklines#1{{\begin{multline*}#1\end{multline*}}} % replace every \cr by \\ 